\documentclass[oneside,english]{amsart}
\pdfoutput=1
\usepackage[T1]{fontenc}
\usepackage[latin9]{inputenc}
\usepackage{verbatim}
\usepackage{float}
\usepackage{amsthm}
\usepackage{amstext}
\usepackage{amssymb}

\makeatletter
\theoremstyle{plain}
\theoremstyle{plain}
\newtheorem{thm}{Theorem}
 \theoremstyle{definition}
  \newtheorem{example}[thm]{Example}
  \theoremstyle{definition}
  \newtheorem{defn}[thm]{Definition}
  \theoremstyle{plain}
  \newtheorem{prop}[thm]{Proposition}
  \theoremstyle{remark}
  \newtheorem{claim}[thm]{Claim}
  \theoremstyle{plain}
  \newtheorem{cor}[thm]{Corollary}

\usepackage{bm}
\usepackage{graphics}
\usepackage{color}
\usepackage{epsfig}
\newcommand{\Hilb}{\mathrm{Hilb}}
\newcommand{\Aut}{\mathrm{Aut}}
\newcommand{\init}{\mathrm{in}}

\makeatother

\usepackage{babel}

\begin{document}

\title{Monomization of Power Ideals and Parking functions}

\author{Craig Desjardins}

\date{18 February 2010}
\begin{abstract}
In this note, we find a monomization of a certain power ideal associated
to a directed graph. This power ideal has been studied in several
settings. The combinatorial method described here extends earlier work of other,
and will work on several other types of power ideals, as will appear in later work.
\end{abstract}
\maketitle

\section{Introduction}

Ideals generated by powers of linear forms have arisen in several
areas recently. They appear in work on linear diophantine equations
and discrete splines \cite{DahmenMicchelli1988}, zonotopla algebra
\cite{holtz-2007}\cite{deconcini-2006}, zonotopal Cox rings \cite{sturmfels-2008},
ideals of fat points \cite{Emsalem19951080}, and other areas. Of
particular interest is the computation of the dimensions of quotients
by these ideals, as well as their hilbert series. In many cases these
computations have been connected to the computation of other statistics
which are more germane to the problem in which they appear (eg.,\cite{ArdilaPostnikov},\cite{sturmfels-2008}). 

In this note we demonstrate a fast algorithm for computing the hilbert
series in an important special case. In particular, the computation
is reduced to the far simpler problem of determining the hilbert series
of certain monomial ideal. In this sense, we have {}``monomialized''
the original ideal. This process also introduces a new notion of parking
function, extending the phenomena which have appeared earlier in the
literature.

\section{Parking Functions and Monomization}

We first recall some definitions and facts about parking functions
and $G$-parking functions. A \emph{parking function} of length $n$
is a sequence of non-negative integers $\left(a_{1},\ldots,a_{n}\right)$
such that for $i$ \[
\#\left\{ j|a_{j}<i\right\} \geq i\]
Note that we allow $a_{i}=0$; the definition is often written for
positive integers instead of non-negative, in which case the `$<$'
above is replaced with `$\leq$'. The definition is equivalent to
requiring that the increasing rearrangement $b_{1}\leq\ldots\leq b_{n}$
has the property $b_{i}<i$. The origin of the term parking function
comes from the following interpretation. Suppose $n$ cars arrive
at a linear parking lot with parking spaces labelled between $0$
and $n-1$. Each car, $i$, has a preferred parking space, $f\left(i\right)$.
The cars arrive in order and drive to their preferred spot. If it
is already taken, then they drive until they reach the next empty
spot and take that one. A preference function $f$ is a parking function
if (and only if) every car gets a parking spot without having to turn
around. 

A $G$-parking function is a generalization of parking functions,
appearing first in \cite{PostnikovShapiro2004}, and later in \cite{Benson20101340,Chebikin200531}.
Let $G$ be a digraph on vertices labelled $0,1,\ldots,n$. We will
call $0$ the root of $G$. For every non-empty subset $I\subset\left[n\right]$
and $i\in I$ define $d_{I}\left(i\right)$ to be the number of edges
originating at $i$ and terminating at a vertex \emph{not} in $I$.
A $G$-\emph{parking function} is defined to be a function $f$ assigning
a non-negative integer to the vertices $1,\ldots,n$ such that for
every non-empty subset $I$ there is an $i\in I$ such that $f\left(i\right)<d_{I}\left(i\right)$.
When $G=K_{n+1}$ the complete graph on $n+1$ vertices, the $G$-parking
functions are the same as the parking functions defined above. 

We now give a reinterpretation of the $G$-parking functions which
was in fact the original motivation for their definition. For any
$I\subset\left[n\right]$, we define $D_{I}$ to be the the total
number of edges of $G$ which originate at a vertex in $I$ and terminate
at a vertex outside of $I$. Explicitly, \[
N_{I}=\sum_{i\in I}d_{I}\left(i\right)\]
Now let $k$ be any integer. If $D_{I}+k>0$ for every $I$ as above,
we define a polynomial $p_{I}$ in the ring $\mathbb{C}\left[x_{1},\ldots,x_{n}\right]$
given by \[
p_{I}=\left(\sum_{i\in I}x_{i}\right)^{D_{I}+k}\]
Note that this will always be the case when $k$ is positive. We then
define $\mathcal{I}_{G,k}$ to be the ideal generated by all such
$p_{I}$, and $\mathcal{A}_{G,k}$ to be the quotient $\mathbb{C}\left[x_{1},\ldots,x_{n}\right]/\mathcal{I}_{G,k}$.
Since $\mathcal{I}_{G,k}$ is a homogeneous ideal, $\mathcal{A}_{G,k}$
will have a basis of monomials. Given a monomial basis $B$ of $\mathcal{A}_{G,k}$,
the set of monomials $\mathcal{M}=\mathbb{C}\left[x_{1},\ldots,x_{n}\right]\backslash B$
is an ideal, and $B$ is the basis of \emph{standard monomials} for
$\mathbb{C}\left[x_{1},\ldots,x_{n}\right]/\mathcal{M}$. We call
any such $\mathcal{M}$ a \emph{monomization }of the ideal $\mathcal{I}_{G,k}$.
Our program is to find a monomization for the ideals $\mathcal{I}_{G,k}$
which is natural in some way and easy to compute. Such a theory would
greatly simplify the study of the linear structure of the rings $\mathcal{A}_{G,k}$. 

In the case $k=0$, the picture is especially beautiful. In \cite{PostnikovShapiro2004},
it is shown that the monomials $x^{a}=x_{1}^{a_{1}}\cdot\ldots\cdot x_{n}^{a_{n}}$
where $\left(a_{1},\ldots,a_{n}\right)$ is a $G$-parking function
give a monomial basis for $\mathcal{A}_{G,0}$. Translating this to
the ideal $\mathcal{J}_{G,0}$ of monomials $x^{a}$ where $a$ is
\emph{not }a $G$-parking function, we find the following description.
For every non-empty $I\subset\left[n\right]$, let \[
m_{I}=x_{i_{1}}^{d_{I}\left(i_{1}\right)}\cdot\ldots\cdot x_{i_{r}}^{d_{I}\left(i_{r}\right)}\]
Then, $\mathcal{J}_{G,0}=\left\langle m_{I}\right\rangle _{I\subset\left[n\right]}$. 

There are several features of this monomization we would like to emulate.
Firstly, a set of generators for $\mathcal{J}_{G,0}$ is relatively
easy to compute from $G$. Secondly, both the generators of $\mathcal{I}_{G,0}$
and $\mathcal{J}_{G,0}$ are indexed by the same set, and in particular
are the same size. There is a third very useful feature of $\mathcal{J}_{G,0}$,
which we would like to emulate, but won't be able to. Clearly the
group $H=\Aut_{*}\left(G\right)$ of basepoint preserving automorphisms
of $G$ acts on $\mathbb{C}\left[x_{1},\ldots,x_{n}\right]$, and
since it preserves the $p_{I}$, $H$ also acts on $\mathcal{A}_{G,k}$.
Since the definition of $G$-parking function is invariant under the
action of $H$, we see that this basis is $H$-invariant. We can thus
use the combinatorial structure of the $G$-parking functions to understand
not just the linear structure of $\mathcal{A}_{G,0}$, but also its
structure as an $H$-representation. For example, this shows that
the multiplicity of the trivial representation of $\mathfrak{S}_{n}$
on $\mathcal{A}_{K_{n+1},0}$ is equal to $\frac{1}{n+1}{2n \choose n}$,
the $n^{th}$ Catalan number. 

For $k\neq0$, we cannot in general find an $\Aut_{*}\left(G\right)$
invariant basis of monomials, even when $G$ is a complete graph. 
\begin{example}
Let $G=K_{3}$, the triangle. Then $\mathcal{I}_{G,0}=\left(x^{2},y^{2},\left(x+y\right)^{2}\right)$.
The parking functions are $\left\{ \left(0,1\right),\left(1,0\right),\left(0,0\right)\right\} $,
and indeed\[
\mathcal{A}_{G,0}=\mathbb{C}\oplus\mathbb{C}x\oplus\mathbb{C}y\]
Now for $k=1$, $\mathcal{I}_{G,1}=\left(x^{3},y^{3},\left(x+y\right)^{3}\right)$.
The monomials which are non-zero in $\mathcal{A}_{G,1}$ are $1,x,y,x^{2},xy,y^{2},x^{2}y,$
and $xy^{2}$. As we can easily verify, however, the Hilbert series
of $\mathcal{A}_{G,1}$ is $1+2t+3t^{2}+t^{3}$. In particular, any
monomial basis must contain $1,x,y,x^{2},xy,$ and $y^{2}$, and must
can contain exactly one of $x^{2}y$ or $xy^{2}$. Thus there is no
way to choose an $\mathfrak{S}_{3}$-invariant basis of monomials. 
\end{example}
We now address the case $k=1$. Similar to the $k=0$ case, we define
a monomial $m_{I}$ for any non-empty subset $I\subset\left[n\right]$.
Namely, let \[
\nu_{I}\left(i\right)=\begin{cases}
d_{I}\left(i\right)+1 & \text{\ensuremath{i\in I} is minimal}\\
d_{I}\left(i\right) & \text{\ensuremath{i\in I} is not minimal}\\
0 & i\notin I\end{cases}\]
and define \[
m_{I}=\prod_{i=1}^{n}x_{i}^{\nu_{I}\left(i\right)}\]
Thus these monomials, as in the case $k=0$, are as close to the center
of the Newton polytope of $p_{I}$ as possible. Let $\mathcal{J}_{G,1}=\left\langle m_{I}\right\rangle $.
The main result of this note is that $\mathcal{J}_{G,1}$ is a monomization
of $\mathcal{I}_{G,1}$. 
\begin{thm}
\label{thm:Monomization of k=00003D1}The standard monomial basis
of $\mathcal{J}_{G,1}$ give a monomial basis for $\mathcal{I}_{G,1}$ 
\end{thm}
We can use this to give a more combinatorial defintion of what we
will call a $\left(G,1\right)$-parking function. In the case of the
complete graph, we believe this definition appeared first in \cite{holtz-2007}.
\begin{defn}
For a graph $G$ on the vertex set $\left\{ 0,1,\ldots,n\right\} $,
a $\left(G,1\right)$-parking function is a function $f:\left[n\right]\to\mathbb{N}$
such that for any $I\subset\left[n\right]$, \[
f\left(i\right)<\begin{cases}
\#\left\{ \text{edges from \ensuremath{i} out of \ensuremath{I}}\right\} +1 & \text{ if \ensuremath{i} is minimal}\\
\#\left\{ \text{edges from \ensuremath{i}out of \ensuremath{I}}\right\}  & \text{otherwise}\end{cases}\]
\end{defn}
\begin{example}
Let's take the following example.

\begin{figure}[H]
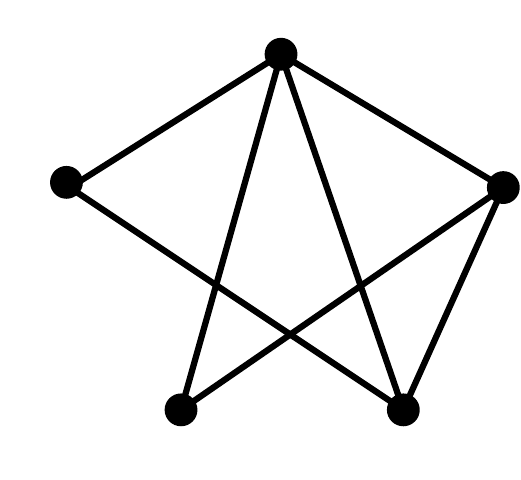

\caption{}

\end{figure}
From the graph above we obtain the power ideal\begin{alignat*}{1}
\mathcal{I}_{G,1}=(x_{1}^{3},x_{2}^{3}, & x_{3}^{4},x_{4}^{4},\left(x_{1}+x_{2}\right)^{5},\left(x_{1}+x_{3}\right)^{4},\left(x_{1}+x_{4}\right)^{6},\left(x_{2}+x_{3}\right)^{6},\\
 & \left(x_{2}+x_{4}\right)^{4},\left(x_{3}+x_{4}\right)^{5},\left(x_{1}+x_{2}+x_{3}\right)^{6},\left(x_{1}+x_{2}+x_{4}\right)^{6},\\
 & \left(x_{1}+x_{3}+x_{4}\right)^{5},\left(x_{2}+x_{3}+x_{4}\right)^{5},\left(x_{1}+x_{2}+x_{3}+x_{4}\right)^{5})\end{alignat*}
and the monomization\begin{alignat*}{1}
\mathcal{J}_{G,1}=(x_{1}^{3},x_{2}^{3},x_{3}^{4} & ,x_{4}^{4},x_{1}^{3}x_{2}^{2},x_{1}^{2}x_{3}^{2},x_{1}^{3}x_{4}^{3},x_{2}^{3}x_{3}^{3},x_{2}^{2}x_{4}^{2},\\
 & x_{3}^{3}x_{4}^{2},x_{1}^{2}x_{2}^{2}x_{3}^{2},x_{1}^{3}x_{2}x_{4}^{2},x_{1}^{2}x_{3}x_{4}^{2},x_{2}^{2}x_{3}^{2}x_{4},x_{1}^{2}x_{2}x_{3}x_{4})\end{alignat*}
From this, it is fairly easy to compute that the dimension of $\mathcal{B}_{G,1}$,
and therefore $\mathcal{A}_{G,1}$, is equal to $82$. Notice, it
is fairly easy to reduce the number of generators in $\mathcal{J}_{G,1}$.
For example, clearly since we have $x_{1}^{3}$, we don't need $x_{1}^{3}x_{2}^{3}$,
$x_{1}^{3}x_{4}^{3}$, or $x_{1}^{3}x_{2}x_{4}^{2}$. Continuing in
this way we can reduce to a minimal set of $10$ generators. 
\end{example}
It is worth mentioning the initial ideal of $\mathcal{I}_{G,k}$ and
its differences from the ideals $\mathcal{J}_{G,0}$ and $\mathcal{J}_{G,1}$.
Given a term order, i.e. a linear ordering on the monomials, we can
form the initial ideal $\init\left(\mathcal{I}_{G,k}\right)$, the
set of leading terms of every element of $\mathcal{I}_{G,k}$. A Gröbner
basis for $\mathcal{I}_{G,k}$ is a set of generators $\left\{ f_{s}\right\} _{s\in S}$
such that the leading terms of the $f_{s}$ generate $\init\left(\mathcal{I}_{G,k}\right)$.
The theory of Gröbner bases is very rich and general, and provides
algorithms for studying a surprising amount of structure of any ideal
in a polynomial ring. However, determining a Gröbner basis for a general
ideal is a potentially time intensive procedure. Additionally, Gröbner
bases can be much larger than a given set of generators for an ideal. 

The ideals $\mathcal{J}_{G,0}$ and $\mathcal{J}_{G,1}$ are almost
never initial ideals of $\mathcal{I}_{G,0}$ and $\mathcal{I}_{G,1}$.
Indeed, initial ideals are generated by vertices of the newton polytopes
of a Gröbner basis for the ideal whereas our ideals are generated
by choosing monomials near the center of the newton polytope. The
tradeoff is that the ideals $\mathcal{J}_{G,k}$ are much easier to
compute, relying only on valence data of the associated graph. Furthermore,
the size of their sets of generators is strictly controlled ($2^{n}-1$),
and typically smaller than that of a Gröbner basis. 
\begin{example}
If we let $G=K_{5}$, the $15$ monomial generators of $\mathcal{J}_{G,1}$
are all non-redundant. However, this is still an improvement over
the $26$ elements of a Gröbner basis for for $\mathcal{I}_{G,1}$,
and are much more difficult to compute. 
\end{example}

\section{Monotone Monomial Ideals}

For the remainer of the paper, let $\mathcal{B}_{G,1}=\mathbb{C}\left[x_{1},\ldots,x_{n}\right]/\mathcal{J}_{G,1}$.
Our proof will proceed as follows. We first demonstrate the standard
monomial basis of $\mathcal{B}_{G,1}$ spans $\mathcal{A}_{G,1}$.
We then show that $\dim\left(\mathcal{B}_{G,1}\right)=\dim\left(\mathcal{A}_{G,1}\right)$,
from which we can conclude the result. In order to show the dimensions
are equal, we show that the dimension of $\mathcal{B}_{G,1}$ is equal
to the number of forests on $G$, and use the equivalent result for
$\mathcal{A}_{G,1}$, due to Ardila and Postnikov.
\begin{thm}
\cite{ArdilaPostnikov} The dimension of the algebra $\mathcal{A}_{G,1}$
is equal to the number of forests on $G$. 
\end{thm}
Both of the results will follow from the fact that the set $\left\{ m_{I}\right\} $
is a monotone monomial family, in the language of \cite{PostnikovShapiro2004}.
We will only need the simplest part of this machinery, which we recall
now. Let $\left\{ m_{I}\right\} $ be any collection of monomials
in $\mathbb{C}\left[x_{1},\ldots,x_{n}\right]$. Let $m_{J\backslash I}$
be the monomial formed from $m_{J}$ by removing all $x_{i}$ with
$i\in I$, and let $\bar{I}=\left[n\right]\backslash I$. Then the
collection $\left\{ m_{I}\right\} $ is a monotone monomial family
if $m_{I\backslash I}=1$ and if $I\subset J$, then $m_{J\backslash\bar{I}}$
divides $m_{I}$. 

It is routine to check that the monomials defined above for $k=1$
are a monotone monomial family. Indeed, the condition $m_{I\backslash I}=1$,
which simply states that $m_{I}$ contains $x_{i}$ only if $i\in I$,
is satisfied by definition. To check the second condition, we examine
the degree of $x_{i}$ for $i\in I$ in $m_{I}$ and $m_{J}$. Since
$J\supset I$, the number of edges originating at vertex $i$ and
terminating outside $J$ is smaller than those terminating outside
$I$, ie. $d_{J}\left(i\right)\leq d_{I}\left(i\right)$. So there
are two cases. If $i$ is the smallest element of $I$, then either
its degree goes from $d_{I}\left(i\right)+1$ to $d_{J}\left(i\right)+1$
(in the case that $i$ is still the smalelst element of $J$), or
it goes from $d_{I}\left(i\right)+1$ to $d_{J}\left(i\right)$. In
both cases the degree drops. If $i$ is not the smallest element of
$I$, then it can't be the smallest element of $J$, so the degree
goes from $d_{I}\left(i\right)$ to $d_{J}\left(i\right)$. Therefore,
in any case $\deg_{x_{i}}\left(m_{I}\right)\geq\deg_{x_{i}}\left(m_{J}\right)$. 

Using this we can conclude from \cite[Theorem 3.1]{PostnikovShapiro2004}
our first claim that the standard monomial basis of $\mathcal{B}_{G,1}$
spans $\mathcal{A}_{G,1}$. To investigate the dimension of the of
space $\mathcal{B}_{G,1}$, we first use \cite[Prop. 8.4]{PostnikovShapiro2004}
to find an expression for the dimension as an alternating sum. Let
$\nu_{I}\left(i\right)=\deg_{x_{i}}\left(m_{I}\right)$.
\begin{prop}
The dimension of $\mathcal{B}_{G,1}$ is equal to the alternating
sum \begin{equation}
\sum_{I_{1}\subsetneq\ldots\subsetneq I_{k}}\left(-1\right)^{k}\prod_{i\in I_{1}}\left(\nu\left(i\right)-\nu_{I_{1}}\left(i\right)\right)\times\ldots\times\prod_{i\in I_{k}\backslash I_{k-1}}\left(\nu\left(i\right)-\nu_{I_{k}}\left(i\right)\right)\times\prod_{i\notin I_{k}}\nu\left(i\right)\label{eq:Alternating Sum Dimension}\end{equation}

where we include the empty chain of subsets where $k=0$. 
\end{prop}
We give the following interpretation to the alternating sum. For a
given chain $I_{1}\subsetneq\ldots\subsetneq I_{k}$ the product counts
the number of directed subgraphs $H$ of $G$ with the following properties
\begin{enumerate}
\item There is at most one edge originating at each $i\in\left[n\right]$,
and there is no edge originating at $0$. 
\item If $i\in I_{j}$ for some $j$, then the edge originating at $i$
must end in $I_{j}$ as well.
\item If $i\in I_{j}$ is the minimal element of $I_{j}$, then $i$ has
an edge originating at it. 
\end{enumerate}
Let us note some properties of these subgraphs. Firstly, any subgraph
of $G$ satisfying the first condition appears in the sum with $k=0$.
Secondly, we can embed the set of forests canonically in this collection
as follows. For any forest $F$ of $G$, orient each edge of $F$
so that each connected component has a unique sink at the minimal
element of that component. (Insert example). Note that each such graph
appears as above with $k=0$, and \emph{only }with $k=0$; case analysis
here. Thus the alternating sum counts each forest exactly once. 

We claim that every other subgraph is cancelled out in the sum, so
that the alternating sum is precisely equal to the number of subforests
of $G$. To show this, we now construct an involution on the set of
pairs $\left(H,I_{1}\subsetneq\ldots\subsetneq I_{k}\right)$ . The
involution will only act on the chain of subsets, that is, it will
leave $H$ fixed. A pair will be fixed by the involution if and only
if it corresponds to $\left(H,\emptyset\right)$ with $H$ a canonically
oriented forest. and it will take a chain of length $k$ either to
chain of length $k-1$ or length $k+1$. Since there are no fixed
points, this will show that any nonforest $H$ is cancelled out in
the alternating sum.

The involution will only use some subset of the vertices of $H$ which
we call \emph{special}. We use the following algorithm to label the
vertices of $H$ special and non-special.
\begin{itemize}
\item Let $v$ be the smallest unlabelled vertex. 
\item If $v$ has an edge originating at it, label it and all remaining
unlabelled vertices special and stop. Otherwise, label $v$ non-special
as well as any vertex such that the chain of edges originating from
it terminates at $v$.
\item Return to the first step.
\end{itemize}
Notice that $0$ will always be chosen first, when none of the vertices
are labelled. We also have the following claim.

\begin{claim}
Each connected component is either composed entirely of non-special
vertices or special vertices. Those labelled non-special are trees
rooted at their minimal element.\end{claim}
\begin{proof}
Suppose that $i$ is non-special. Then $i$ must lie on a directed
path towards a terminal vertex, and in particular the path originating
at $i$ does not contain a circuit. This follows because the only
way a vertex can be labelled non-special is in step $2$ of the algorithm,
and only as part of a path which terminates. Therefore, if $i$ is
non-special, then its connected component must be a tree. 

If $i$ is part of a tree $T$ and non-special, we claim the tree
is rooted, ie. it has a unique sink. More specifically, it is rooted
at the end of the path originating from $i$. If this is the case,
then every vertex of the tree was labelled non-special in the same
step that $i$ was. To see that it has a uniqe sink, let $w$ be a
sink in $T$. There is a unique undirected path $\left(w,v_{1},v_{2},\ldots,v_{k},i\right)$
from $w$ to $i$. The edge from $w$ to $v_{1}$ must be oriented
towards $w$, because $w$ is a sink. Because each vertex can have
at most one out-edge, this is the unique out-edge from $v_{1}$. Therefore
the edge from $v_{2}$ to $v_{1}$ must be oriented towards $v_{1}$. 

Continuing in this way, we conclude that the entire path is oriented
from $i$ to $w$. Therefore $w$ is the vertex at the end of the
path originating from $i$, and consequently is unique. The only thing
left to see is that $w$ is the minimal vertex of the tree. This is
easy, though, since otherwise $w$ would not have been chosen in step
$1$ of the algorithm. 
\end{proof}
Note that the converse of the second part of the claim is false. It
is perfectly feasible for a subtree of $H$ to be oriented towards
its minimal vertex and still be labelled special. We do however get
the following corollary. 
\begin{cor}
\label{cor:non-special forest}$H$ is a canonically oriented forest
if and only if the algorithm labels every vertex non-special.\end{cor}
\begin{proof}
One direction is clear from the claim: If all the vertices are labelled
non-special, then every connected component is a rooted tree oriented
toward the minimal vertex, which is a canonically oriented forest.
For the converse, suppose $H$ is a canonically oriented forest, but
some vertex is marked special. Then at some point in the algorithm,
the least unlabelled vertex $v$ has an out-edge. The path coming
out of $v$ must terminate at a vertex smaller than $v$ since the
forest is canonically oriented, but then this vertex must have been
labelled non-special. This, in turn, would imply that $v$ is labelled
non-special, and we get a contradiction. 
\end{proof}

\begin{cor}
If $i\in I_{j}$, then $i$ is special. \end{cor}
\begin{proof}
Let $w_{j}$ be the minimal vertex in $I_{j}$. Clearly $w_{j}$ must
be special, because if it were non-special then by claim blah blah
it would lie on a tree oriented towards its minimal vertex. However
this isn't possible since there is an edge originating at $w_{j}$,
and the entire path from $w_{j}$ must lie within $I_{j}$ by definition.
However, if $i$ is non-special, then the path originating at $i$
must terminate at a non-special vertex $v$, and that vertex must
lie in $I_{j}$. Because $v$ must also lie in $I_{j}$ it must be
greater than $w_{j}$, but then $w_{j}$ would have been chosen in
step $1$ of the above algorithm before $v$ was, and $v$ wouldn't
have been marked non-special. This is a contradiction.
\end{proof}
We now define the involution $\kappa$. Let $\mathcal{S}$ be the
set of special vertices in $H$. Then \[
\kappa\left(\left(H,I_{1}\subsetneq I_{2}\subsetneq\ldots\subsetneq I_{k}\right)\right)=\begin{cases}
\left(H,I_{1}\subsetneq I_{2}\subsetneq\ldots\subsetneq I_{k}\subsetneq\mathcal{S}\right) & \text{ if }\mathcal{S}\neq I_{k}\\
\left(H,I_{1}\subsetneq I_{2}\subsetneq\ldots\subsetneq I_{k-1}\right) & \text{ if }\mathcal{S}=I_{k}\end{cases}\]
If the chain of subsets is empty, and there are no special vertices,
then $\kappa$ does nothing. By Corollary \ref{cor:non-special forest}
this means that $H$ is a canonically oriented tree. Otherwise, the
length of the chain of subsets is changed by $\kappa$. Therefore,
the only fixed points are $\left(H,\emptyset\right)$ where $H$ is
a canonically oriented tree. Applying $\kappa$ to formula \ref{eq:Alternating Sum Dimension}
we get 
\begin{thm}
The dimension of $\mathcal{A}$ is equal to the number of forests
on $G$.
\end{thm}
This completes the proof of \ref{thm:Monomization of k=00003D1}.
The proof actually gives a little more. Since the standard monomials
of $\mathcal{B}_{G,1}$ span $\mathcal{A}_{G,1}$, we get the inequality
\begin{equation}
\Hilb\left(\mathcal{B},t\right)\leq\Hilb\left(\mathcal{A},t\right)\label{eq:Hilbert Inequality}\end{equation}
We have just seen that $\Hilb\left(\mathcal{A},1\right)=\dim_{k}\left(\mathcal{A}\right)$
is equal to the number of forests in $G$, which implies the corollary.
\begin{cor}
$\Hilb\left(\mathcal{B},t\right)=\Hilb\left(\mathcal{A},t\right)$
\end{cor}
From \cite{ArdilaPostnikov} we also have a combinatorial interpretation
for the coeffients in $\Hilb\left(\mathcal{A},t\right)$. If \[
\Hilb\left(\mathcal{A},t\right)=\sum c_{n}t^{n}\]
then $c_{n}$ is equal to the number of forests $F$ on $G$ with
external activity equal to $\left|G\right|-\left|F\right|-n$. This
generalizes the results of \cite{PostnikovShapiroShapiro1998,PostnikovShapiroShapiro1999}.
It would be nice to find a bijection between the monomials of fixed
degree and the forests of fixed external activity.

\section{Further Work}

Some interesting questions remain about these ideals. As noted above,
the monomial ideals $\mathcal{J}_{G,1}$ are \emph{not} in general
the initial ideals of the ideals $\mathcal{I}$. However, since their
Hilbert Series are equal, they each correspond to a point on the same
Hilbert Scheme. We can then ask whether or not they lie on the same
irreducible component of the Hilbert Scheme. If not, we can ask how
many irreducible components away they are from each other. Note that
this number exists, since the Hilbert scheme is connected. 

The original motivation for the study of these rings comes from differential
geometry. Let $G=SL_{n}\left(\mathbb{C}\right)$ and $B$ the subgroup
of matrices fixing a given flag in $\mathbb{C}^{n}$, then $G/B$
is the complex flag variety, parametrizing complete flags in $\mathbb{C}^{n}$.
There are $n$ tautological line bundles over $G/B$, assigning to
each point $p\in G/B$ the quotient of the $k$-plane by the $\left(k-1\right)$-plane
of the flag corresponding to $p$. Fixing a hermitian structure on
$G/B$ gives us a unique connection on each of these line bundles,
and therefore curvature forms $\omega_{1},\ldots,\omega_{n}$. Note
that the deRham cohomology classes of these $2$-forms (after normalizing)
are precisely the Chern classes of the corresponding line bundles.
For this reason, we can think of the ring $\mathbb{C}\left[\omega_{1},\ldots,\omega_{n}\right]$,
sitting inside the ring of $C^{\infty}$ invariant forms on $G/B$,
as an extension of the intersection theory obtained from $H^{*}\left(G/B,\mathbb{C}\right)$.
In \cite{PostnikovShapiroShapiro1999} it is shown that this ring
is isomorphic to the ring $\mathcal{A}_{K_{n+1},1}$, and therefore
has a basis given by the monomials given above. It would be interesting
to determine if other rings of Chern forms on homogeneous spaces have
similar presentations, and to determine monomizations of their ideals.

\bibliographystyle{plain}
\bibliography{biblio}

\end{document}